\documentclass[12pt, twoside]{article}
\usepackage{amsmath,amsthm,amssymb}
\usepackage{times}
\usepackage{enumerate}

\def\titlerunning#1{\gdef\titrun{#1}}
\makeatletter
\def\author#1{\gdef\autrun{\def\and{\unskip, }#1}\gdef\@author{#1}}
\def\address#1{{\def\and{\\\hspace*{18pt}}\renewcommand{\thefootnote}{}%
\footnote {#1}}%
\markboth{\autrun}{\titrun}}
\makeatother
\def\email#1{e-mail: #1}
\def\subjclass#1{{\renewcommand{\thefootnote}{}%
\footnote{\emph{Mathematics Subject Classification (2020):} #1}}}

\numberwithin{equation}{section}

\begin{document}

\baselineskip=17pt

\titlerunning{Development of sieve of Eratosthenes and sieve of Sundaram's proof}

\title{Development of sieve of Eratosthenes and sieve of Sundaram's proof}

\author{Ahmed Diab}

\date{}

\maketitle

\address{
\email{ahmedhamdy1401753@sci.asu.edu.eg}
}


\begin{abstract}
We make two algorithms that generate all prime numbers up to a given limit, they are a development of sieve of Eratosthenes algorithm, we use two formulas to achieve this development, where all the multiples of prime number 2 are eliminated in the first formula, and all the multiples of prime numbers 2 and 3 are eliminated in the second formula.

Using the first algorithm we  proof sieve of Sundaram's algorithm, then we improve it to be more efficient prime generating algorithm.

We will show the difference in performance between all the algorithms we will make and sieve of Eratosthenes algorithm in terms of run time.

\end{abstract}

\section{Introduction}
Sieve of Eratosthenes algorithm (SOE) is a well known algorithm that generates all prime numbers up to a given limit $N$, it contains some rules that control the elimination process of the multiples of prime numbers see \cite{10.2307/41186370,10.2307/3607073,horsley1772kosigma,o2009genuine,sorenson1990introduction}, two of these rules are
\begin{enumerate}
\item \begin{equation}
N_f=R^2
\end{equation} 
Where $N_f$ is the first uncommon multiple of prime number $R$ with prime numbers smaller than $R$, which will be the first multiple to be eliminated.
\item
\begin{equation}
R_l=\left\lfloor\sqrt{N}\right\rfloor
\end{equation}
Where $R_l$ is the last prime number whose multiples will be eliminated.
\end{enumerate} 
Using equations (3.2) and (3.4) in reference  \cite{diab2021sequence} and SOE algorithm we will develop two algorithms that generate all prime numbers up to a given limit $N$.

We use the advantage of equation (3.2) in reference \cite{diab2021sequence} which is 
\begin{equation}
2n+1=p
\end{equation}
Where all multiples of 2 are eliminated,
so we will set our rules to eliminate the orders of the multiples of prime numbers greater that 2 from $n$, then the remaining values of $n$ will be entered to equation (1.3) to get its corresponding prime numbers.\\

We will use the first produced algorithm as a proof to Sieve of Sundaram's algorithm \cite{aiyar1934sundaram} ,then we will develop it using the property of the exact match of the mathematical concepts between these two algorithms.\\

Then we will use equation (3.4) in reference \cite{diab2021sequence} to make the second development to SOE algorithm, this is the formula
\begin{equation}
2\left\lfloor\frac{3n_1+1}{2}\right\rfloor+1=p
\end{equation}
Where all multiples of 2 and 3 are eliminated, we will set some rules to eliminate the orders of the multiples of prime numbers greater than 3 from $n_1$ .

\section{First development of SOE (D1SOE)}
\label{sec:1st algorithm }

In this developed algorithm we eliminate the corresponding $n$ of the multiples of  prime numbers greater than 2 from $n$ - in equation (1.3) - up to a certain limit $n_m$, where
\begin{equation}
n_m=\left\lfloor\frac{N-1}{2}\right\rfloor
\end{equation} 
We take the floor to the right hand-side of equation (2.1) because we are interested in the integer value of $n_m$.

To obtain the order of the first uncommon multiple of a prime number with prime numbers smaller we substitute by ($N_f=2m+1$) and ($R=2z+1$) in equation (1.1)

\begin{center}
$2m+1=(2z+1)^2$
\end{center}
\begin{center}
$2m+1=4z^2+4z+1$
\end{center}
\begin{equation}
m=2(z^2+z)
\end{equation}

$m$: Is the order of the first multiple we will eliminate.

$z$: Is the order of the prime number whose multiples we want to eliminate.

In this case the multiples of prime numbers are equidistant numbers, consequently their orders too see \cite{diab2021sequence}, in order to eliminate the multiples of prime numbers we need to know the difference between the orders of any two successive multiples, we can get the multiples of any number by multiplying this number by the integer numbers, but in this case we will multiply it by the odd numbers because the multiples of 2 are already eliminated, so for $p$ prime number, the multiples of $p$ are $((1*p), (3*p), (5*p), (7*p), ... )$.

using first two successive multiples which are $(1p)$ and $(3p)$ then substituting by each one in the inverse of equation (1.3) then subtracting the corresponding $n$ of $(1p)$ from the corresponding $n$ of $(3p)$ we obtain
\begin{center}
$\frac{3p-1}{2}-\frac{p-1}{2}=d$
\end{center} 
\begin{equation}
\therefore p=d
\end{equation}    
$d$: The difference between the orders of any two successive multiples of a prime number $p$.

We will eliminate all the numbers starting from $m$ with $d$ steps to $n_m$ which are the orders of the multiples of prime number $p$.

The order of the last prime number its multiples should be eliminated will be derived from equation (1.2) -without the floor- by substituting by 
\begin{itemize}
\item $R_l=2k+1$
\item $N=2q+1$
\end{itemize}

equation (1.2) will be 
\begin{equation}
k=\frac{\sqrt{2q+1}}{2}-\frac{1}{2}
\end{equation}

We will take the floor to the right hand-side of equation (2.4) because we are interested in the integer value of $k$

\begin{equation}
k=\left\lfloor\frac{\sqrt{2q+1}}{2}-\frac{1}{2}\right\rfloor
\end{equation}

$k$: The order of the number whose first prime number smaller than it is the last prime number whose multiples we want to eliminate.

$q$:  The order of the limit which we want to generate all prime numbers up to.

After eliminating all the orders of prime numbers starting from 3 to the first prime number smaller than $k$ from $n$ we substitute by the remaining values of $n$ in equation (1.3) to obtain all prime numbers in the interval $\left[ 3,N \right]$.\\

Following a c++ code representation for D1SOE algorithm.

\rule{\textwidth}{1pt}
\begin{verbatim}
#include <iostream>
#include <cstdlib>
#include<cmath>
using namespace std;
int D1SOE(int n_m) {
	int n = 0, m = 0, ArrayLim = 0;
	cout << "\n2 ";
	// adding 1 to n_m because the array starts from order 0
	ArrayLim = n_m + 1;
	bool* array = new bool[ArrayLim];
	// Initializing the array with false values 
	for (int i = 0; i < ArrayLim; i++)
		array[i] = false;
	// The mean elimination theorem 
	for (n = 1; n <= int((sqrt(((2 * n_m) + 1) / 4.0)) - 0.5); n++)
	{if (array[n] != 0)
			continue;
		for (m = (2 * ((n * n) + n)); m <= n_m; m += ((2 * n) + 1))
		{array[m] = true;}}
	// Printing for loop
	for (int n = 1; n < ArrayLim; n++)
		if (!array[n])
			cout << (2 * n) + 1 << " ";
	return 0;}
// Driver program to test above 
int main() {
	int n_m, N = 0;
	cout << "\n Enter limit : ";
	cin >> N;
	n_m = ((N - 1) / 2);
	D1SOE(n_m);
	return 0;}
\end{verbatim}
\rule{\textwidth}{1pt}
\section{Proof and development of sieve of Sundaram (DSOS)}
\label{sec:SOS proof}
let $i,j\in\mathbb{N},\ 1\le i\le j$
\begin{equation}
i+j+2ij=y
\end{equation} 
Sieve of Sundaram's algorithm \cite{aiyar1934sundaram} eliminates all the values of $y$ from the integer numbers up to a limit $l$ where $y\le l$.

We can write equation (3.1) as follows
\begin{equation}
i+(2i+1)j=y
\end{equation}
Let 
\begin{equation}
u=2i+1
\end{equation}


Condition $i\le j$ will control the starting point of $y$, which is the value of $y$ when $i=j$.

substituting by $i=j$ in equation (3.2)
\begin{center}
$i+(2i+1)i=y$
\end{center}
\begin{equation}
2(i^2+i)=y
\end{equation}   
Equation (3.4) is the same as equation (2.2) which also determines the starting point of the elimination process.

The end point of the elimination process is $y\le l$ where $l$ by the definition of sieve of Sundaram is the integer value of 
\begin{center}
$l=\frac{N-1}{2}$
\end{center}
which could be written as 
\begin{equation}
l=\left\lfloor\frac{N-1}{2}\right\rfloor
\end{equation}

$N$: Is the limit which we want to generate all the prime numbers up to.

Equation (3.5) is the same as equation (2.1) which determines the end point of the elimination process, this condition will control the last value of $j$.\\
Starting from ($2(i^2+i)$) to $l$ with $u$ steps represented by the term ($u*j$) in equation (3.2) this algorithm is nearly the same as D1SOE algorithm.\\
Substitutes by the remaining integer values in an equation the same as equation (1.3) prime numbers are obtained.

The differences between D1SOE algorithm and Sundaram's algorithm is
\begin{enumerate}
\item In equation (3.1) $i$ takes all integer values from 1 to $l$, in D1SOE it takes only the values of the orders of prime numbers.
\item In Sundaram's algorithm the limit of $i$ is defined by equation (3.5), where D1SOE uses equation (2.5).

Using equation (2.5) this is the limit of $i$ should be used
\begin{equation}
k'=\left\lfloor\frac{\sqrt{2q'+1}}{2}-1\right\rfloor
\end{equation}
where\\ 
$k'$: Is the order of -the corresponding $n$ of- the limit of $i$.

$q'$: Is the order of -the corresponding $n$ of- $N$.
\end{enumerate}
 Changing these conditions in Sundaram's algorithm it will be exactly the same as D1SOE algorithm.\\
 
Following a c++ code representation for DSOS algorithm.

 \rule{\textwidth}{1pt}
\begin{verbatim}
#include <iostream> 
#include <cstdlib>
using namespace std;
int DSOS(int n)
{   int i = 0, ArrayLim = 0, j = 0, n_m = 0;
    n_m = (n - 1) / 2;
       cout << 2 << " ";
    // adding 1 to n_m because the array starts from order 0
    ArrayLim = n_m + 1;
    bool* array = new bool[ArrayLim];
    // Initializing the array with false values
    for (i = 0; i < ArrayLim; i++)
        array[i] = false;
    // The mean elimination theorem         
    for (i = 1; i<=int((sqrt(((2 * n_m) + 1) / 4.0)) - 0.5); i++)
    {   if (array[i] != 0)
            continue;
        for (j = i; (i + j + 2 * i * j) <= n_m; j++)
            array[i + j + 2 * i * j] = true; }
      // Printing for loop
      for (i = 1; i <= n_m; i++)
        if (array[i] == false)
            cout << 2 * i + 1 << " ";
    return 0;}
// Driver program to test above 
int main(void)
{   int n = 0;
    cout << "\n Enter limit : ";
    cin >> n;
    DSOS(n);  
    return 0;}
\end{verbatim}
\rule{\textwidth}{1pt}
 
\section{Second development of SOE algorithm (D2SOE)}
\label{sec:2nd algorithm }
In this algorithm we eliminate the orders of the multiples of prime numbers greater than 2 and 3 from $n_1$ in equation (1.4).

Using ISEF from reference \cite{diab2021sequence} the inverse of equation (1.4) will be
\begin{center}
$n_1=\left\lceil\frac{((\frac{p-1}{2})*2)-1}{3}\right\rceil$
\end{center}
\begin{equation}
n_1=\left\lceil\frac{p-2}{3}\right\rceil
\end{equation}

The limit of $n_1$ is

\begin{equation}
n_{1m}=\left\lceil\frac{N-2}{3}\right\rceil
\end{equation}
$N$: The limit which we want to generate all the prime numbers up to.\\
$n_{1m}$: The order of $N$.\\

The order of the first multiple we want to eliminate will be calculated by substituting by  
\begin{itemize}
\item $\left\lfloor\frac{3g+1}{2}\right\rfloor=m$
\item $\left\lfloor\frac{3b+1}{2}\right\rfloor=z$
\end{itemize}
in equation (2.2)
\begin{center}
$\left\lfloor\frac{3g+1}{2}\right\rfloor=2(\left\lfloor\frac{3b+1}{2}\right\rfloor^2+\left\lfloor\frac{3b+1}{2}\right\rfloor)$
\end{center}

Using ISEF from reference \cite{diab2021sequence}
\begin{equation}
g=\left\lceil\frac{4(\left\lfloor\frac{3b+1}{2}\right\rfloor^2+\left\lfloor\frac{3b+1}{2}\right\rfloor)-1}{3}\right\rceil
\end{equation}
$g$: The order of the first multiple will be eliminated.\\
$b$: The order of the prime number whose multiples we want to eliminate.\\

The multiples of the prime numbers in equation (1.4) are divided into two sequences see \cite{diab2021sequence}, and $g$ can equal any number from these two sequences.

The multiples of a prime number $p$ in equation (1.4) are ($p$, $5p$, $7p$, $11p$,...) which are the integer numbers without the multiples of 2 and 3 multiplied by prime number $p$.
The first sequence is ($p$, $7p$, $13p$, ...), and the second sequence is ($5p$, $11p$, $17p$, ...).

We need to know the following information to be able to eliminate the multiples of the prime numbers in this case  
\begin{enumerate}
\item The difference between the orders of any two successive multiples in any sequence of these two sequences, let it be $f_1$.
\begin{equation}
f_1=\left\lceil\frac{7p-2}{3}\right\rceil-\left\lceil\frac{p-2}{3}\right\rceil
\end{equation}
Using computer programming the following relation is true for $p$ equals the integers from 1 to $10^{10}$
\begin{equation}
2p=\left\lceil\frac{7p-2}{3}\right\rceil-\left\lceil\frac{p-2}{3}\right\rceil
\end{equation}
We assume that this relation is true for $p$ equals all integers upto infinity. 
 
\item The difference between the orders of the first (or the second or the third, ...etc) numbers in both sequence, let it be $f_2$.
\begin{equation}
f_2=\left\lceil\frac{5p-2}{3}\right\rceil-\left\lceil\frac{p-2}{3}\right\rceil
\end{equation}
\item The difference between the order of the first number in the second sequence and the order of the second number in the first sequence, let it be $f_3$.
\begin{equation}
f_3=\left\lceil\frac{7p-2}{3}\right\rceil-\left\lceil\frac{5p-2}{3}\right\rceil
\end{equation}
\item Whether $g$ belongs to the first or the second sequence, the second sequence could be represented by $(3x-1)$ where $x$ is an integer, so by adding 1 to $g$ and then dividing it by 3 we can know that it belongs to the second sequence if the result is an integer, if the result is not an integer then it belongs to the first sequence, while the first sequence could be represented by $(3x+1)$.
\end{enumerate}

After calculating $g$ we will determine whether $g$ belongs to the first or the second sequence, Then we will eliminate this sequence by adding $f_1$ to $g$ then we will eliminate the result then we will repeat adding $f_1$ to the last result and eliminate the resultant number, our limit will be $n_{1m}$.

After eliminating the sequence which $g$ belongs to we will eliminate the other sequence by the same method but starting from ($g+f_3$) if $g$ belongs to the second sequence, and from (($g+f_2$) or ($g-f_3$)) if $g$ belongs to the first sequence.\\

The order of the last prime number its multiples should be eliminated will be derived from equation (2.4) by substituting by 
\begin{itemize}
\item $\left\lfloor\frac{3h+1}{2}\right\rfloor=k$
\item $\left\lfloor\frac{3e+1}{2}\right\rfloor=q$
\end{itemize}

Equation (2.4) will be
 \begin{equation}
 \left\lfloor\frac{3h+1}{2}\right\rfloor=\frac{\sqrt{2\left\lfloor\frac{3e+1}{2}\right\rfloor+1}}{2}-\frac{1}{2}
 \end{equation}
 using ISEF 
 \begin{equation}
h=\left\lceil\frac{\sqrt{2\left\lfloor\frac{3e+1}{2}\right\rfloor+1}-2}{3}\right\rceil
 \end{equation}
 where\\ 
 $h$: The order -corresponding $n_1$- of the number whose first prime number smaller than it is the last prime number whose multiples we want to eliminate.\\
$e$:  The order -corresponding $n_1$- of the limit which we want to generate all prime numbers up to.

After eliminating the orders of the multiples of prime numbers from $n_1$ we will substitute by the remaining $n_1$ in equation (1.4) to obtain prime numbers.\\

Following a c++ code representation for D2SOE algorithm.

 \rule{\textwidth}{1pt}
\begin{verbatim}
#include <iostream>
#include <cstdlib>
using namespace std;
int D2SOE(int n1_m) {
	int n1 = 0, g = 0, z = 0, p = 0,f1=0,f3=0;
	cout << "\n2 3 ";
	bool* array = new bool[n1_m];
	// Initialising the D2SOE array with false values 
	for (int i = 0; i < n1_m; i++)
		array[i] = false;
	// The mean elimination theorem 
	for (n1=1;n1<=ceil((sqrt(2*floor((3*n1_m+1)/2.0)+1)-2)/3.0);n1++)
	{ if (array[n1] != 0)
			continue;
		z = ((3 * n1 + 1) / 2.0);
		p = ((2 * z) + 1);
		f1 = (ceil(((7*p)-2) / 3.0) - ceil((p - 2) / 3.0));
		f3 = (ceil(((7*p)-2) / 3.0) - ceil(((5 * p) - 2) / 3.0));
		for (g = ceil(((4*((z*z) + z)) - 1) / 3.0); g < n1_m;g+=f1)
		{ array[g] = true; }
		if ((p +1) % 3 == 0) {
			for (g = ceil(((4*((z*z) + z)) - 1) / 3.0) + f3; g < n1_m;g+=f1)
			{ array[g] = true; } }
		else {
			for (g = ceil(((4*((z*z) + z)) - 1) / 3.0) - f3; g<n1_m;g+=f1)
			{ array[g] = true; }}}
	// printing for loop
	 for (int n1 = 1; n1 < n1_m; n1++)
		if (!array[n1]) {
			z = (((3 * n1) + 1) / 2.0);
			cout << (2 * z) + 1 << " "; }
	return 0; }
// driver program to test above
int main() {
	int n1_m, N = 0;
	cout << "\n Enter limit : ";
	cin >> N;
	n1_m = ceil((N - 2) / 3.0);
	D2SOE(n1_m);
	return 0; }
\end{verbatim}
\rule{\textwidth}{1pt}

\section{Implementation results}
\label{sec:impl.results}

To know in general the speed difference between these algorithms we use the same coding concept for all algorithms and the same device under the same circumstances, the results are listed in the next table, where this table shows the run time in microseconds for each algorithm required to generate all the prime numbers from 1 to the limit $N$. 

\begin{center}
\begin{tabular}{ c | c | c | c | c }
\hline $N$ & SOE & D1SOE & D2SOE & DSOS \\
\hline
$10^3$ & 10 & 10 & 10 & 9\\
$10^4$ & 92 & 53 & 32 & 54\\
$10^5$ & 978 & 472 & 263 & 450\\
$10^6$ & 12793 & 6021 & 3492 & 5812\\
$10^7$ & 190523 & 88483 & 53848 & 86673\\
$10^8$ & 2160761 & 1045056 & 741409 & 1038911\\
$10^9$ &24837596 & 12143628 & 8471540 & 12008482\\  
$2*10^9$ &..... & 25136075 & 17953830 & 25208467\\
\hline
\end{tabular}
\end{center}

D2SOE is the fastest algorithm then D1SOE then DSOS with almost the same speed as D1SOE, then the last and the slowest one is SOE algorithm.

\section{conclusion}
Using sieve of Eratosthenes algorithm and two of my prime numbers formulas we was able to make two developments to sieve of Eratosthenes algorithm, the first one we used a formula where all the multiples of 2 are eliminated and we derived all the required formulas to accomplish the elimination of the multiples of the prime numbers, and  the generation of all the prime numbers up to a gives limit see \ref{sec:1st algorithm }, and we used this algorithm to proof sieve of Sundaram's algorithm by showing the exact match between these two algorithms, then we used this match to improve sieve of Sundaram see section \ref{sec:SOS proof}.
 
In the second algorithm we used a formula where all the multiples of the prime numbers 2 and 3 are eliminated, and derived the required formulas for the elimination process see section \ref{sec:2nd algorithm }.

Finally we showed the run time differences between our developed algorithms and sieve of Eratosthenes algorithm see section \ref{sec:impl.results}.

\bibliographystyle{plain}
\bibliography{references}

\end{document}